\makeatletter
\@ifundefined{extract@font}
{
\documentstyle[12pt,leqno]{article}

\newcommand{\text}[1]{\mbox{\rm #1}}
\newcommand{\Bbb}[1]{\mbox{\bf #1}}
\newcommand{\operatorname}[1]{\mathop{\rm #1}\nolimits}

\newcommand{\pfit}[1]{{\it #1:\/}}
\newcommand{\qedsymbol}{Q.E.D.}
\newcommand{\eqref}[1]{(\ref{#1})}
\newcommand{\dotsb}{\cdots}
\newcommand{\binom}[2]{{{#1}\choose{#2}}}

\newenvironment{pf}{\noindent\pfit{Proof} }{\qedsymbol\par\par\medskip\par}
\newenvironment{pf*}{\noindent\pfit }{\qedsymbol\par\par\medskip\par}

\newenvironment{bmatrix}{\left[\begin{array}{ccccc}}{\end{array}\right]}
\title{
Picard-Fuchs equations and mirror maps \\ for hypersurfaces}
\author{David R. Morrison}
\newcommand{\trailer}{\medskip\par\noindent {\sc
Department of Mathematics,
Duke University,
Durham, NC\ \ 27706
}\par\noindent {\it E-mail:}
drm@math.duke.edu}
\date{}
}
{
\documentstyle[12pt]{amsart}
\title[Picard-Fuchs equations and mirror maps]{
Picard-Fuchs equations and mirror maps \\ for hypersurfaces}
\author{David R. Morrison}
\address{
Department of Mathematics \\
Duke University \\
Durham, NC\ \ 27706}
\email{drm@math.duke.edu}
\newcommand{\trailer}[1]{}
}
\makeatother

\newtheorem{lemma}{Lemma}  
\newtheorem{corollary}{Corollary}  
\newtheorem{proposition}{Proposition} 
\newtheorem{theorem}{Theorem} 
\renewcommand{\qedsymbol}{Q.E.D.}

\newcommand{\spn}{\operatorname{span}}
\newcommand{\suchthat}{\ | \ }
\newcommand{\Pic}{\operatorname{Pic}}
\newcommand{\Res}{\operatorname{Res}}
\newcommand{\lcm}{\operatorname{lcm}}
\newcommand{\mugroup}[1]{\mbox{\boldmath ${\mu}_{#1}$}}

\begin{document}

\maketitle

\setlength{\unitlength}{1in}
\hfill
\begin{picture}(1.1,0)
\put(0,2.2){\makebox{\parbox{1.1in}{{\sc duk-m-91-14} \\
                                         October, 1991}}}
\end{picture}

\begin{abstract}
We describe a strategy for computing Yukawa couplings and
the mirror map, based on the Picard-Fuchs equation.
(Our strategy is a variant of the method used by
Candelas, de la Ossa, Green, and Parkes \cite{pair} in the
case of quintic hypersurfaces.)
We then explain a technique of Griffiths \cite{griffiths} which can
be used to compute
the Picard-Fuchs equations of hypersurfaces.
Finally,  we carry out the computation for
four specific examples (including quintic hypersurfaces, previously
done by Candelas et al.~\cite{pair}).  This yields predictions
for the number of rational curves of various degrees on certain
hypersurfaces in weighted projective spaces.
Some of these predictions have been confirmed by
classical techniques in algebraic geometry.
\end{abstract}

\section*{Introduction}

The phenomenon of mirror symmetry dramatically caught the attention
of mathematicians with the recent work of
P.~Candelas, X.~C.~de la Ossa, P.~S.~Green, and L.~Parkes \cite{pair}.
Starting with a particular pair of ``mirror manifolds'',
 calculating certain period integrals,
interpreting the results as Yukawa couplings, and then re-interpreting
those results in light of the ``mirror manifold'' phenomenon,
Candelas et al.\ were
able to give predictions for the numbers of rational curves of various
degrees on the general quintic threefold.
In fact, algebraic geometers have had a difficult time
verifying these predictions, but all
successful attempts to calculate the numbers of curves
have eventually confirmed the predictions.

What is so striking about this work
 is that the calculation which predicts the
numbers of rational curves on quintic threefolds
 is in reality a calculation about the
variation of Hodge structure on a completely {\em different\/} family of
Calabi-Yau threefolds.  An asymptotic expansion is made of a function
which comes from that variation, and the coefficients in the
expansion are then used to predict numbers of rational curves.

In \cite{guide}, we interpreted the calculation of Candelas et
al.~\cite{pair}
in terms of variation of Hodge structure.  Here we take a more down
to earth approach, and work directly with period integrals and their
properties.  (This is perhaps closer in spirit to the original paper.)
We have found a way to modify the computational strategy employed in
\cite{pair}.  Our modified method computes a bit less (there are
two unknown ``constants of integration''), but it
is easier to actually carry out the computation.  We in
fact carry it out in three new examples.  This leads to new predictions
about numbers of rational curves on certain Calabi-Yau threefolds.

Our strategy for computing Yukawa couplings is based on the Picard-Fuchs
equation for the periods of a one-parameter family of algebraic varieties.
We explain in sections 1 and 2 how this equation can be used to compute
Yukawa couplings and
the mirror map for a family of Calabi-Yau threefolds with
 $h^{2,1}=1$.  We then go on in section 3 to
review a method of Griffiths \cite{griffiths} for calculating
Picard-Fuchs equations of hypersurfaces.  Related ideas have also been
introduced into the physics literature in
\cite{blok-var,cad-ferr,ferrara,lsw}.

In sections 4 and 5, we carry out the computation in four examples,
including the quintic hypersurface.  The resulting
predictions about numbers of
rational curves are discussed in section 6.

\section{The Picard-Fuchs equation and monodromy}

Let $\bar\pi: \overline{\cal X}\to \overline C$
 be a  family of $n$-dimensional
projective algebraic varieties,
parameterized by
 a compact Riemann surface  $\overline C$.  Let $C \subset \overline C$
be an open subset such that the induced family $\pi: {\cal X}\to C$
has smooth fibers.  If we choose
topological $n$-cycles $\gamma_0,\dots,\gamma_{r-1}$ which give a basis
for the $n^{\text{th}}$
homology of one particular
fiber $X_0$, and choose  a holomorphic $n$-form $\omega$ on $X_0$,
then the {\em periods} of $\omega$ are the integrals
\[ \int_{\gamma_0}\omega,\dots,\int_{\gamma_{r-1}}\omega. \]

Since the fibration $\pi: {\cal X}\to C$ is differentiably
locally trivial, a local trivialization can be used to extend
the cycles $\gamma_i$ from $X_0$ to cycles $\gamma_i(z)$ on $X_z$
which depend on
$z$, where $z$ is a local coordinate on $C$.
The holomorphic $n$-form $\omega$ can also be extended to a family of
$n$-forms $\omega(z)$ which depend on the parameter $z$.  If this is
done in an algebraic way, then $\omega(z)$ extends to a meromorphic
family of
$n$-forms (i.e. poles are allowed) over the entire space $\overline{\cal X}$.

The cycles $\gamma_i(z)$ determine homology classes which
are locally constant in $z$.
However, an attempt to extend these cycles
 globally will typically lead to monodromy:
for each closed path in $C$,
there will be some linear map $T$ represented by a matrix $T_{ij}$
such that transporting $\gamma_i$ along the path produces at the end
a cycle homologous to
$\sum T_{ij}\gamma_j$.
The same phenomenon will hold for the periods:  for a globally
defined meromorphic family of $n$-forms $\omega(z)$, the local periods
 $\int_{\gamma_i(z)}\omega(z)$ extend by analytic continuation to
 multiple-valued functions of $z$, transforming according to the same
monodromy transformations
$T$ as do the homology classes of the cycles.

The periods $\int_{\gamma(z)}\omega(z)$
 satisfy an ordinary differential equation called
the {\em Picard-Fuchs equation\/} of $\omega$.
 The existence of this equation can be explained as follows.
Choose a local coordinate $z$
on some open set $U\subset C$, and consider the vector
\[ v_j(z) := [\frac{d^j}{dz^j}\int_{\gamma_0(z)}\omega(z),\dots,
\frac{d^j}{dz^j}\int_{\gamma_{r-1}(z)}\omega(z)] \in\Bbb C^r. \]
For generic values of the parameter $z$,
the dimensions
\[
d_j(z):=\dim(\spn\{v_0(z),\dots,v_j(z)\})
\]
must be constant.
Since $d_j(z)\le r$, these spaces cannot continue to grow indefinitely.
There
will thus be a smallest $s$ such that
\[v_s(z)\in\spn\{v_0(z),\dots,v_{s-1}(z)\}\]
(for generic $z$).
We can write
\[ v_s(z) =-\sum_{j=0}^{s-1}C_j(z)v_j(z) \]
with the coefficients $C_j(z)$ depending on $z$.
The {\em Picard-Fuchs equation},
satisfied by all the periods of $\omega(z)$, is then
\begin{equation} \label{picfuc}
 \frac{d^sf}{dz^s}+\sum_{j=0}^{s-1}C_j(z)\frac{d^jf}{dz^j}=0.
\end{equation}
The precise form of the equation depends on both the local coordinate
$z$ on $C$,
and the choice of holomorphic form $\omega(z)$.
Note that the coefficients $C_j(z)$ may acquire singularities at special
values of $z$.

When we approach a point $P$
in $\overline C-C$, the Picard-Fuchs equation has (at worst)
a {\em regular singular point\/} at $P$ \cite{gr-bull,nkatz,deligne}.
If we choose a parameter $z$ which is
centered at $P$ (that is, $z=0$ at $P$), then the coefficients $C_j(z)$
in the Picard-Fuchs equation typically will have poles at $z=0$.
However, if we multiply the Picard-Fuchs operator
\begin{equation} \label{PFop}
 \frac{d^s}{dz^s}+\sum_{j=0}^{s-1}C_j(z)\frac{d^j}{dz^j}
\end{equation}
by $z^s$ and rewrite the result in the form
\begin{equation} \label{logform}
 (z\frac{d}{dz})^s+\sum_{j=0}^{s-1}B_j(z)(z\frac{d}{dz})^j
\end{equation}
then the new coefficients $B_j(z)$ are holomorphic functions of $z$.
(This is one of several equivalent definitions of ``regular singular point''.)
We call eq.~\eqref{logform} the {\em logarithmic form\/} of the Picard-Fuchs
operator.

The structure of ordinary differential equations with regular singular points
is a classical topic in differential equations:  a convenient reference is
\cite{codlev}.  We can rewrite eq.~\eqref{picfuc} as a system of first-order
equations, using the logarithmic form eq.~\eqref{logform}, as follows:
let
\begin{equation} \label{Az}
 A(z) = \begin{bmatrix}
    0    &    1    &        &        &            \\
         &    0    &   1    &        &            \\
         &         & \ddots & \ddots &            \\
         &         &        &   0    &    1       \\
 -B_0(z) & -B_1(z) & \dots  & \dots & -B_{s-1}(z)
\end{bmatrix}.
\end{equation}
Then solutions $f(z)$ to the equation eq.~\eqref{picfuc} are equivalent
to solution
vectors
\[w(z) = \begin{bmatrix} f(z) \\ z\frac d{dz}f(z) \\ \vdots \\
(z\frac d{dz})^{s-1}f(z) \end{bmatrix} \]
of the matrix equation
\begin{equation} \label{matrixeqn}
z\frac d{dz}w(z)=A(z)w(z).
\end{equation}

For a matrix equation such as eq.~\eqref{matrixeqn}, the facts are these
(see \cite{codlev}).
  There is
a constant $s\times s$ matrix $R$ and a $s\times s$ matrix $S(z)$
of (single-valued)
functions of $z$, regular near $z=0$, such that
\[\Phi(z)=S(z) \cdot z^R \]
is a {\em fundamental matrix} for the system.  This means that the
columns of
$\Phi(z)$ are a basis for the space of solutions at each nonsingular
point $z\ne0$.
The multiple-valuedness of the solutions has all been put into $R$, since
\[z^R:=e^{(log z)R} = I+(\log z)R+\frac{(\log z)^2}{2!}R^2+\dotsb\]
is a multiple-valued matrix function of $z$.  The local monodromy on
the solutions
given by analytic continuation along a path winding
 once around $z=0$ in a counterclockwise direction is given
by $e^{2\pi iR}$ (with respect to the basis given by the columns of $\Phi$).
The matrix $R$ is by no means unique.

\begin{theorem}
Suppose that $z\frac d{dz}w(z)=A(z)w(z)$ is a system of ordinary
differential
equations with a regular singular point at $z=0$.  Suppose that
distinct eigenvalues
of $A(0)$ do not differ by integers.  Then there is a fundamental
matrix of
the form
\[\Phi(z)=S(z) \cdot z^{A(0)} \]
and $S(z)$ can be obtained as a power series
\[S(z)=S_0+S_1z+S_2z^2+\dotsb\]
by recursively solving the equation
\[
z\frac d{dz}S(z)+S(z)\cdot A(0)=A(z)\cdot S(z)
\]
for the coefficient matrices $S_j$.
Moreover, any such series solution
converges in a neighborhood of $z=0$.
\end{theorem}
A proof can be found in \cite{codlev}, together with
 methods for treating the case in which eigenvalues of
$A(0)$ {\em do} differ by integers.

We will be particularly interested in systems with
{\em unipotent monodromy}:
by definition, this means that $e^{2\pi iR}$ is a unipotent matrix,
so that
$(e^{2\pi iR}-I)^m\ne0$, $(e^{2\pi iR}-I)^{m+1}=0$ for some $m$
called the
{\em index}.

\begin{corollary}
Suppose that
$ (z\frac{d}{dz})^sf(z)+\sum_{j=0}^{s-1}B_j(z)(z\frac{d}{dz})^jf(z) $
is an ordinary differential equation with a regular singular point at $z=0$.
If $B_j(0)=0$ for all $j$, then the solutions of this equation have
unipotent
monodromy of index $s$.
\end{corollary}
The corollary follows by calculating with
eq.~\eqref{Az}, setting $z=0$ and $B_j(0)=0$
to produce
\[
e^{2\pi iA(0)}= \begin{bmatrix}
1 & 2\pi i & \frac{(2\pi i)^2}{2!} & \dots & \frac{(2\pi i)^{s-1}}{(s-1)!} \\
  & 1 & 2\pi i & \dots & \frac{(2\pi i)^{s-2}}{(s-2)!} \\
 & & \ddots & & \vdots \\
 & & & 1 & 2\pi i  \\
 & & & & 1
\end{bmatrix}.\]

\section{Computing the mirror map}

Recall that a {\em Calabi-Yau manifold} is a compact K\"ahler manifold $X$
of complex dimension $n$ which has trivial canonical bundle,
such that the Hodge numbers $h^{k,0}$ vanish for $0<k<n$.
Thanks to a celebrated theorem of Yau \cite{yau}, every such manifold
admits Ricci-flat K\"ahler metrics.

Suppose now that $\pi: {\cal X}\to C$ is a family of Calabi-Yau threefolds
with $h^{2,1}(X)=1$, which is not a locally constant family.
The third cohomology group $H^3(X)$ has dimension
$r=4$.  It follows that the Picard-Fuchs equation has order at most 4.
(In fact, it is not difficult to show
that it has order exactly 4.)

Let $z$ be a coordinate on $\overline C$ centered at a point
$P\in\overline C-C$.  We say that
{\em $P$ is a point at which the monodromy is maximally
unipotent} if the monodromy is unipotent of index 4.  As we have seen
in the corollary, if $B_j(0)=0$ in the logarithmic form of
the Picard-Fuchs equation, $z=0$ will be such a point.  We will assume
for simplicity that our points of maximally unipotent monodromy have
this form, leaving appropriate modifications for the general case to
the reader.

We review the calculation of the Yukawa coupling, following \cite{pair}.
Let $\omega(z)$ be a family of $n$-forms, and let
\[W_k := \int_{X_z}\omega(z)\wedge\frac{d^k}{dz^k}\omega(z).\]
A fundamental principle from the theory of variation of Hodge structure
(cf.~\cite{transcendental}) implies that
$W_0$, $W_1$, and $W_2$ all vanish.  The {\em Yukawa coupling} is the first
non-vanishing term $W_3$.  Candelas et al.\ show that the Yukawa coupling
$W_3$ satisfies the differential equation
\[ \frac{dW_3(z)}{dz}=-\frac12C_3(z)W_3(z), \]
where $C_3(z)$ is a coefficient in the Picard-Fuchs equation \eqref{picfuc}.

The Yukawa coupling as defined clearly depends on the ``gauge'', that is,
on the choice of holomorphic $3$-form $\omega(z)$.  If fact, if we alter
the gauge by $\omega(z)\mapsto f(z)\omega(z)$, then $W_k$ transforms
as
\[ W_k\mapsto f(z)\sum_{j=0}^k\binom{k}{j}\frac{d^jf(z)}{dz^j}W_{k-j}. \]
Since $W_0=W_1=W_2=0$, the change in the Yukawa coupling $W_3$ is
simply $W_3\mapsto f(z)^2W_3$.

The Yukawa coupling also depends on the choice of coordinate $z$,
and in fact is often denoted by $\kappa_{zzz}$.
If we change coordinates from $z$ to $w$, we must change the differentiation
operator from $d/dz$ to $d/dw$.  The chain rule then
imples that
\[ \kappa_{www} = \left(\frac{dz}{dw}\right)^3\kappa_{zzz}.\]

Candelas et al.~\cite{pair} use physical arguments to set the gauge in this
calculation, and to find an appropriate (multiple-valued) parameter
 $t$ with which to compute.
(The  associated
differentiation operator $d/dt$ is single-valued.)
What will be important for us are
the following observations about their results.

The gauge used by Candelas et al.\ determines a family of meromorphic
$n$-forms $\widetilde\omega(z)$ with the property that the period function
\[\int_{\gamma}\widetilde\omega(z)\equiv1\]
for some cycle $\gamma$.  Moreover, the parameter $t$ determined by Candelas
et al.\ is a parameter defined in an angular sector near $z=0$
which has two crucial
properties:
\begin{enumerate}
\item
If we analytically continue along
 a simple loop around $z=0$ in the counterclockwise direction, $t$
becomes $t+1$.
(It will be convenient to also introduce $q=e^{2\pi it}$, which
remains single-valued near $z=0$.)

\item
There are cycles $\gamma_0$ and $\gamma_1$ such that
$\int_{\gamma_0}\omega(z)$
is single valued near $z=0$, and
\[t=\frac{\int_{\gamma_1}\omega(z)}{\int_{\gamma_0}\omega(z)}\]
in an angular sector near $z=0$.
\end{enumerate}

Each period function $\int_{\gamma}\omega(z)$ is a solution to the
Picard-Fuchs equation of the family.  Translating the results of
the previous section into the present context, we obtain the following:

\begin{lemma}
Suppose that $z=0$ is a point of maximally unipotent monodromy
such that $B_j(0)=0$, where $B_j(z)$ are the coefficients in the
logarithmic form of the Picard-Fuchs equation.  Then
\begin{enumerate}
\item
There is a period function for $\omega(z)$,
\[f_0(z):=\int_{\gamma_0}\omega(z)\]
which is single-valued near $z=0$.  This period function is unique up to
multiplication by a constant.
(This implies that the cycle $\gamma_0$ is also unique up to a
constant multiple.)

In particular, the family of meromorphic $n$-forms
\[\widetilde\omega(z):=\frac{\omega(z)}{\int_{\gamma_0}\omega(z)}\]
will have the property that
\[\int_{\gamma}\widetilde\omega(z)\equiv1\]
for some $\gamma$, and it is the unique such family up to constant multiple.

\item
Fixing a choice of period function $f_0(z)$ as in part (1), there is a
period function
\[f_1(z):=\int_{\gamma_1}\omega(z)\]
such that $\varphi(z):=f_1(z)/f_0(z)$ transforms as
\[\varphi(z)\mapsto\varphi(z)+1\]
upon transport around $z=0$ in the counterclockwise direction.
The ratio $\varphi(z)$ is unique up to the addition of a constant.
\end{enumerate}
\end{lemma}

This, then, is our alternate strategy for computing the Yukawa coupling:
we find solutions of the Picard-Fuchs equation which have the properties
specified in the lemma, and we use those to fix the gauge and specify the
natural parameter, up to two unknown constants of integration.

\section{Picard-Fuchs equations for hypersurfaces}

We now review a method of Griffiths \cite{griffiths} for describing
the cohomology of a hypersurface, which can be used to determine
the Picard-Fuchs equation of a one-parameter family of hypersurfaces.
Calculations of this sort were earlier made by Dwork \cite[Sec. 8]{dwork}.
Griffiths' method was extended to the weighted projective case by
Steenbrink \cite{steen} and Dolgachev \cite{dolg}, who we follow.

We denote a weighted projective $n$-space by $\Bbb P^{(k_0,\dots,k_n)}$,
where $k_0,\dots,k_n$ are the weights of the variables $x_0,\dots,x_n$.
Weighted homogeneous polynomials can be identified with the aid of the
Euler vector field
\[
\theta=\sum k_jx_j\frac{\partial}{\partial x_j}
\]
which has the property that $\theta P=(\deg P)\cdot P$ for any weighted
homogeneous polynomial $P$.  Contracting the volume from on $\Bbb C^{n+1}$
with $\theta$ produces the fundamental weighted homogeneous differential
form (of ``weight'' $k:=\sum k_j$)
\[
\Omega:=\sum_{j=0}^n(-1)^jk_jx_j\,
dx_0\wedge\dots\wedge\widehat{dx_j}\wedge\dots\wedge dx_n.
\]
Rational differentials of degree $n$ on $\Bbb P^{(k_0,\dots,k_n)}$
can be described as expressions $P\Omega/Q$, where $P$ and $Q$ are
weighted homogeneous polynomials with $\deg P+k=\deg Q$.

Suppose that $Q$ is a weighted homogeneous polynomial defining a quasismooth
hypersurface ${\cal Q}\subset \Bbb P^{(k_0,\dots,k_n)}$.
(That is, $Q=0$ defines a hypersurface in $\Bbb C^{n+1}$ which is smooth
away from the origin.)  The middle cohomology of ${\cal Q}$ is then
described by means of differential forms with poles (of all orders)
along ${\cal Q}$.  Each such form $P\Omega/Q^\ell$ is made into a
cohomology class by a ``residue'' construction:  for an $(n-1)$-cycle
$\gamma$ on ${\cal Q}$, the tube over $\gamma$ (an $S^1$-bundle inside
the (complex) normal bundle of ${\cal Q}$) is an $n$-cycle $\Gamma$
on $\Bbb P^{(k_0,\dots,k_n)}$ disjoint from ${\cal Q}$.  We can then
define the residue of $P\Omega/Q^\ell$ by
\[
\int_{\gamma}\Res_{\cal Q}\left(\frac{P\Omega}{Q^\ell}\right)=
\frac1{2\pi i}\int_{\Gamma}\frac{P\Omega}{Q^\ell}.
\]
Since altering $P\Omega/Q^\ell$ by an exact differential does not change
the value
of these integrals, we see that the cohomology of ${\cal Q}$ is
represented by
equivalence classes of rational differential forms $P\Omega/Q^\ell$
modulo exact forms.

Here is Griffiths' ``reduction of pole order'' calculation which shows
how to reduce modulo exact forms in practice.  Let $Q$ and $A_j$
be weighted homogeneous polynomials, with $\deg Q=d$,
$\deg A_j =\ell d +k_j-k$.  Define
\[
\varphi=\frac1{Q^\ell}\sum_{i<j}(k_ix_iA_j-k_jx_jA_i)
dx_0\wedge\dots\wedge\widehat{dx_i}\wedge
\dots\wedge\widehat{dx_j}\wedge\dots\wedge dx_n
\]
and then calculate
\begin{equation} \label{griffiths}
d\varphi=\frac
{\left(\ell\sum A_j\frac{\partial Q}{\partial x_j} -Q\sum\frac{\partial A_j}
{\partial x_j}\right)\Omega}{Q^{\ell+1}}
=\frac{\ell\sum A_j\frac{\partial Q}{\partial x_j}\Omega}{Q^{\ell+1}}
-\frac{\sum\frac{\partial A_j}
{\partial x_j}\Omega}{Q^\ell}.
\end{equation}
Thus, any form whose numerator lies
in the Jacobian ideal
$
 J=(\partial Q/\partial x_0,
\dots,\partial Q/\partial X_n)
$
is equivalent (modulo exact forms)
to a form with smaller pole order.

This idea can be used to calculate Picard-Fuchs equations as follows.
The cycles $\Gamma$ do not change (in homology) when $z$ varies locally.
So we can differentiate under the integral sign
\[
\frac{d^k}{dz^k}\int_{\gamma}\Res_{\cal Q}\left(\frac{P\Omega}{Q^\ell}\right)=
\frac1{2\pi i}\int_{\Gamma}\frac{d^k}{dz^k}\left(\frac{P\Omega}{Q^\ell}\right)
\]
when $Q$ depends on a parameter $z$.  (Note that $\Omega$ is independent
of $z$.)  The Picard-Fuchs operator \eqref{PFop} will have the property
that
\[
\left( \frac{d^s}{dz^s}+\sum_{j=0}^{s-1}C_j(z)\frac{d^j}{dz^j} \right)
\left(\frac {P\Omega}Q\right)=d\varphi
\]
is an exact form.  To find it, take successive $z$-derivatives of the
integrand
$P\Omega/Q$ and use the reduction of order of pole formula \cite{griffiths}
to determine a linear relation among those derivatives, modulo exact forms.

\section{Examples: Picard-Fuchs equations}

We will calculate the Picard-Fuchs equations for
certain one-parameter families of Calabi-Yau threefolds.
Our choice of families is motivated by
 the mirror construction of Greene and Plesser \cite{greene-plesser}.

We choose weights $k_0,\dots,k_4$ with $k_0\ge k_1\ge\dots\ge k_4$
for a weighted projective 4-space
such that $d_j:=k/k_j$ is
an integer, where $k:=\sum k_j$.
We also assume that $\gcd\{k_j\suchthat j\ne j_0\}=1$ for every $j_0$.
These assumptions then imply that $k=\lcm\{d_j\}$.

Consider the pencil of hypersurfaces ${\cal Q}_\psi \subset
\Bbb P^{(k_0,\dots,k_4)}$ defined by $Q(x,\psi)=0$, where
\[
Q(x,\psi) :=  \sum_{j=0}^4 x_j^{d_j} - k\psi\prod_{j=0}^4 x_j  .
\]
This pencil has a natural group of diagonal automorphisms preserving the
holomorphic 3-form.  To define it,
let $\mugroup{m}$ denote the multiplicative group of $m^{\text{th}}$ roots
of unity (considered as a subgroup of $\Bbb C^\times$), and let
\[
G=(\mugroup{d_0}\times\dots\times\mugroup{d_4})/\mugroup{k},
\]
where we embed $\mugroup{k}$ in $\mugroup{d_0}\times\dots\times\mugroup{d_4}$
by
\[
\alpha\mapsto(\alpha^{k_0},\dots,\alpha^{k_4}).
\]
Note that since $\sum k_j=k$, the formula
\[
f(\alpha_0,\dots,\alpha_4)=(\prod \alpha_j)^{-1}
\]
determines a well-defined homomorphism $f:G\to\Bbb C^\times$.
Let $G_0=\ker(f)$.

We can regard $Q(x,\psi)=0$ as defining a hypersurface
${\cal Q}\subset\Bbb P^{(k_0,\dots,k_4)}\times\Bbb C$.
The group $G$ acts on $\Bbb P^{(k_0,\dots,k_4)}\times\Bbb C$ by
\[
(x_0,\dots,x_4;\psi)\mapsto(\alpha_0x_0,\dots,\alpha_4x_4;f(\alpha)\psi)
\]
for $\alpha=(\alpha_0,\dots,\alpha_4)\in G$.  The polynomial $Q(x,\psi)$
is invariant under this action.  Thus, the
action preserves ${\cal Q}$,
and maps ${\cal Q}_\psi$ isomorphically to ${\cal Q}_{f(\alpha)\psi}$.
It follows that the group $G_0$ acts on ${\cal Q}_\psi$ by automorphisms,
and that the induced action of
$G/G_0\cong\mugroup{k}$ establishes isomorphisms between
${\cal Q}_\psi/G_0$ and ${\cal Q}_{\lambda\psi}/G_0$ for
$\lambda\in\mugroup{k}$.

The quotient
space ${\cal Q}_\psi/G$ has only canonical singularities.
By a theorem of Markushevich \cite[Prop.~4]{markushevich} and
Roan \cite[Prop.~2]{roan1}, these singularities can be
resolved to give a Calabi-Yau manifold ${\cal W}_\psi$.
There are choices to be made in this resolution
process; we do not specify a choice.
By another theorem of Roan \cite[Lemma 4]{roan2},
 any two resolutions differ by a sequence of flops.

Note that the differential form $\Omega$ from the previous section
transforms as $\Omega\mapsto(\prod \alpha_j)\Omega$ under the action
of $\alpha\in G$.  Thus, the rational differential
\[
\omega_1=\frac{\psi\Omega}{Q(x,\psi)}
\]
is invariant under the action of $G$; we define
$\omega(\psi)=\Res_{{\cal Q}_\psi}(\omega_1)$.

Since the holomorphic 3-forms $\omega(\psi)$
on ${\cal Q}_\psi$ are invariant on $G_0$,
they induce holomorphic 3-forms on ${\cal W}_\psi$.  Moreover, the
homology group $H_3({\cal W}_\psi)$ contains the $G_0$-invariant part
$H_3({\cal Q}_\psi)^{G_0}$ of the homology of ${\cal Q}_\psi$.  If we know
that the dimensions of these spaces agree, then they will coincide
(at least for homology with coefficients in a field).  In this case,
the periods of ${\cal W}_\psi$ can actually be computed as periods of
the holomorphic form $\omega(\psi)$
on ${\cal Q}_\psi$, over $G_0$-invariant cycles.
Thanks to the isomorphisms between ${\cal Q}_\psi$ and
${\cal Q}_{\lambda\psi}$
for $\lambda\in\mugroup{k}$ and the invariance of the rational differential
$\omega_1$ under $G$, these periods will be invariant under
$\psi\mapsto\lambda\psi$.  In particular, they will be functions of
$z=\psi^{-k}$ alone.

It is likely that the resolutions ${\cal W}_\psi$
of ${\cal Q}_\psi/G_0$ could be chosen so
that the action of $G/G_0$ would lift to isomorphisms between
${\cal W}_\psi$ and ${\cal W}_{\lambda\psi}$.  (We verified this in the case
of quintic hypersurfaces in \cite{guide}.)  In this case, there would
be an actual family of Calabi-Yau threefolds for which $z$ served as
a parameter.  It may be that such resolutions could be constructed by
finding an appropriate partial resolution of ${\cal Q}/G$.  However,
we do not need the existence of
this family to describe the computation of the Yukawa
coupling.

\begin{table}[t]
\renewcommand{\arraystretch}{1.2}
\begin{center}
\begin{tabular}{|c|c|c|} \hline
$k$ & $(k_0,\dots,k_4)$ & $Q(x,\psi)$ \\ \hline
5 & $(1,1,1,1,1)$ & $x_0^5+x_1^5+x_2^5+x_3^5+x_4^5-5\psi x_0x_1x_2x_3x_4$ \\
6 & $(2,1,1,1,1)$ & $x_0^3+x_1^6+x_2^6+x_3^6+x_4^6-6\psi x_0x_1x_2x_3x_4$ \\
8 & $(4,1,1,1,1)$ & $x_0^2+x_1^8+x_2^8+x_3^8+x_4^8-8\psi x_0x_1x_2x_3x_4$ \\
10 & $(5,2,1,1,1)$ & $x_0^2+x_1^5+x_2^{10}+x_3^{10}+x_4^{10}
-10\psi x_0x_1x_2x_3x_4$ \\
\hline \end{tabular} \end{center}
\caption{The hypersurfaces.}
\label{tab1}
\end{table}

We will carry out the computation in four specific examples.  These
come from the lists of
 Candelas, Lynker and  Schimmrigk
 \cite{cls}; they found that there are
exactly four types of hypersurface in weighted projective four-space
which are Calabi-Yau threefolds with Picard number one.
The weights of the space are given in the second column of table~\ref{tab1}.
For each of
those cases, Greene and Plesser's mirror construction \cite{greene-plesser}
yields the family ${\cal W}_\psi$ which we have described above.
And Roan's formula \cite{roanpf} for the Betti numbers verifies that $b_3$
is indeed 4 (with $h^{2,1}=1$).
The remaining columns in table~\ref{tab1} show the value of $k$, and
give the equation $Q(x,\psi)$ explicitly.

We describe the $G_0$-invariant cohomology by means of the rational
differential forms
\[
\omega_\ell:=
\frac{
(-1)^{\ell-1}(\ell-1)!\,\psi^\ell(\prod x_i^{\ell-1})\Omega}{Q(x,\psi)^\ell}.
\]
These are chosen because of the evident $G$-invariance in the numerator;
the coefficients were adjusted so that the formula
\begin{equation} \label{derivs}
-\frac1k\psi\frac d{d\psi}\omega_\ell=-\frac{\ell}k\omega_\ell
+\omega_{\ell+1}
\end{equation}
would not be overly burdened with constants.  We compute with the
differential
operator $-\frac1k\psi\frac d{d\psi}$ because it coincides with
$z\frac d{dz}$.

A basis for the $G_0$-invariant cohomology is then given by the residues of
$\omega_1$, $\omega_2$, $\omega_3$, $\omega_4$.  To compute the
Picard-Fuchs equation, we must find an expression for $\omega_5$
as a linear combination of $\omega_1,\dots,\omega_4$ modulo exact forms.
That expression, combined with \eqref{derivs}, will then yield the
desired differential equation.

We carried out this calculation using the Gr\"obner basis
algorithm \cite{buchberger},
modifying an implementation written in {\sc maple} by Yunliang Yu
(cf.~\cite{yu}).  We first
calculated a Gr\"obner basis for the Jacobian ideal
$
J=(\partial Q/\partial x_0,\dots,\partial Q/\partial x_4)
$,
working in the ring $\Bbb C(\psi)[x_0,\dots,x_4]$ of polynomials whose
coefficients are rational functions of $\psi$.  The reduction of pole
order was then achieved step by step as follows:  given a form $\eta_\ell$,
the residue of a global form with a pole of order $\ell$, we
used the Gr\"obner basis to reduce the numerators of
both $\eta_\ell$ and $\omega_\ell$ to standard form.  We could thus
determine a coefficient $\varepsilon_\ell\in\Bbb C(\psi)$
such that the numerator of
$\eta_\ell-\varepsilon_\ell\omega_\ell$ lies in $J$.
Another application of Gr\"obner basis reduction produced explicit
coefficients
\[
\eta_\ell-\varepsilon_\ell\omega_\ell=\sum A_{\ell j}\frac
{\partial Q}{\partial x_j}.
\]
Then the
Griffiths formula
\eqref{griffiths} determines forms $\varphi_\ell$
and $\eta_{\ell-1}$ such that
\[
\eta_\ell-\varepsilon_\ell\omega_\ell=d\varphi_\ell+\eta_{\ell-1},
\]
and $\eta_{\ell-1}$ has a pole of order $\ell-1$.

Beginning with $\eta_5=\omega_5$ and applying this procedure several times,
one finds
\[
\omega_5=\varepsilon_1\omega_1+\dots+\varepsilon_4\omega_4 +d\varphi.
\]
The results of this computation for our four examples
are summarized in table 2.  The
coefficients $\varepsilon_\ell$ are in fact functions of $z=\psi^{-k}$
(as expected from our earlier discussion), and have been displayed as
such.

\begin{table}[t]
\renewcommand{\arraystretch}{1}
\begin{center}
\begin{tabular}{|c|cccc|} \hline
\rule{0pt}{14pt} $k$ & $\varepsilon_1$ & $\varepsilon_2$ & $\varepsilon_3$ &
$\varepsilon_4$ \\[4pt]
\hline
 & & & & \\
5 & $\displaystyle\frac{1}{625(z-1)}$ &
$\displaystyle\frac{-3}{25(z-1)}$ & $\displaystyle\frac{1}{(z-1)}$ &
$\displaystyle\frac{-2}{(z-1)}$ \\
 & & & & \\
6 & $\displaystyle\frac{1}{324(z-4)}$ & $\displaystyle\frac{-5}{18(z-4)}$ &
$\displaystyle\frac{-(z-50)}{18(z-4)}$
& $\displaystyle\frac{-(z+20)}{3(z-4)}$ \\
 & & & & \\
8 & $\displaystyle\frac{1}{16(z-256)}$
& $\displaystyle\frac{-15(z+256)}{512(z-256)}$ &
$\displaystyle\frac{-5(3z-1280)}{64(z-256)}$
& $\displaystyle\frac{-(3z+1280)}{4(z-256)}$ \\
 & & & & \\
10 & $\displaystyle\frac{5}{4(z-12500)}$
& $\displaystyle\frac{-(7z+37500)}{200(z-12500)}$ &
$\displaystyle\frac{-(7z-62500)}{20(z-12500)}$
& $\displaystyle\frac{-(z+12500)}{(z-12500)}$ \\
 & & & & \\
\hline \end{tabular} \end{center}
\caption{The results of the Gr\"obner basis calculation.}
\label{tab2}
\end{table}

The differential equation for $[\omega_1,\dots,\omega_4]$ determined
by this procedure has the form
\[
z\frac d{dz}
\begin{bmatrix} \omega_1 \\ \omega_2 \\ \omega_3 \\ \omega_4 \end{bmatrix}
=
\begin{bmatrix}
-\frac1k & 1 & 0 & 0 \\
0 & -\frac2k & 1 & 0 \\
0 & 0 & -\frac3k & 1 \\
\varepsilon_1 & \varepsilon_2 & \varepsilon_3 & \varepsilon_4-\frac4k
\end{bmatrix}
\begin{bmatrix} \omega_1 \\ \omega_2 \\ \omega_3 \\ \omega_4 \end{bmatrix}.
\]
To calculate the Picard-Fuchs equation, we must change basis via
\[
{\renewcommand{\arraystretch}{1.5}
\begin{bmatrix} \omega_1 \\ z\frac d{dz}\omega_1 \\ (z\frac d{dz})^2\omega_1 \\
(z\frac d{dz})^3\omega_1 \end{bmatrix}
=
\begin{bmatrix}
1 & 0 & 0 & 0 \\
-\frac1k & 1 & 0 & 0 \\
\frac1{k^2} & -\frac3k & 1 & 0 \\
-\frac1{k^3} & \frac7{k^2} & -\frac6k & 1
\end{bmatrix}
\begin{bmatrix} \omega_1 \\ \omega_2 \\ \omega_3 \\ \omega_4 \end{bmatrix}.
}
\]
This determines an equation in the form \eqref{Az}, with
\begin{equation} \label{Beqn}
{\renewcommand{\arraystretch}{1.5}
\begin{array}{rrr}
B_0(z) & = & -\varepsilon_1(z) - \frac1k\varepsilon_2(z) -
\frac2{k^2}\varepsilon_3(z)
-\frac6{k^3}\varepsilon_4(z)+\frac{24}{k^4} \\
B_1(z) & = & -\varepsilon_2(z)-\frac3k\varepsilon_3(z)-
\frac{11}{k^2}\varepsilon_4(z)
+\frac{50}{k^3} \\
B_2(z) & = & -\varepsilon_3(z)-\frac6k\varepsilon_4(z)+\frac{35}{k^2} \\
B_3(z) & = & -\varepsilon_4(z)+\frac{10}k.
\end{array}
}
\end{equation}
As can be directly verified in each of our cases, $B_j(0)=0$.
It
follows that the monodromy at $z=0$ is maximally unipotent.
(In the case of quintics ($k=5$), this had been shown in \cite{pair};
cf.~\cite{guide}.)

\section{Examples:  Mirror maps}

We next compute the mirror maps for our four examples, based on their
Picard-Fuchs equations.
Expanding eqs.~\eqref{PFop} and \eqref{logform}, one finds that the
coefficient $C_3(z)$ coincides with $(6+B_3(z))/z$.  Moreover, in
our four examples, a straightforward computation based on eq.~\eqref{Beqn}
and table~\ref{tab2} shows that
 $B_3(z)=2z/(z-\lambda)$,
where $\lambda=1$, $4$, $256$, $12500$ when $k=5$, $6$, $8$, $10$,
respectively.  Thus,
\[
C_3(z)=\frac{6+B_3(z)}{z}=\frac6z+\frac2{z-\lambda}.
\]
The Yukawa coupling $\kappa_{zzz}$ in the gauge $\omega(z)$ is
therefore given by
a function $W_3(z)$ which satisfies the differential equation
\[
\frac{dW_3(z)}{dz}=\left(\frac{-3}{z}+\frac{-1}{z-\lambda}\right)W_3(z).
\]
Thus, in the gauge $\omega(z)$ we have
\[
\kappa_{zzz}=\frac{c_1}{(2\pi i)^3z^3(z-\lambda)}.
\]
Here $c_1/(2\pi i)^3$ is the first ``constant of integration'':
we have introduced a factor of $(2\pi i)^3$ in order to simplify a
later formula.

In order to determine the natural gauge, we must find a solution
$f_0(z)$ of the
Picard-Fuchs equation which is regular near $z=0$.  Using the corresponding
vector $w_0(z)$ of which $f_0(z)$ is the first component,
we want a solution to the vector equation
\begin{equation} \label{w0eqn}
z\frac d{dz}w_0(z)=A(z)w_0(z)
\end{equation}
which is regular near $z=0$.  ($A(z)$ is given by
eqs.~\eqref{Az}, \eqref{Beqn},
and table~\ref{tab2}.)
This can be found
using power-series techniques, and there is a solution with $f_0(0)\ne0$
in each of our four cases.  We  normalize so that $f_0(0)=1$;
alternatively, we could have absorbed the leading term of $f_0(z)$ into the
constant of integration $c_1$.

As a result, the gauge-fixed value of $\kappa_{zzz}$ takes the form
\[
\kappa_{zzz}=\frac{c_1}{(2\pi i)^3z^3(z-\lambda)(f_0(z))^2},
\]
where the constant $c_1$ has yet to be determined.

We now search for the good parameter $t$.  We should locate a second
solution $f_1(z)$, or its corresponding vector $w_1(z)$, which is
multiple-valued and has the correct monodromy properties.
The monodromy will be such that if we introduce
\[v(z):=2\pi iw_1(z)-(\log z)w_0(z)\]
 and its first
component
\[g(z):=2\pi if_1(z)-(\log z)f_0(z),\]
then $v(z)$ will be single-valued and regular
near $z=0$.  It is easy to calculate that the matrix equation
satisfied by $v(z)$ is
\begin{equation} \label{veqn}
z\frac d{dz}v(z)=A(z)v(z)-w_0(z).
\end{equation}
Solutions to this equation
can be found by power-series techniques.  We normalize the solution
so that $g(0)=0$.  The parameter $t$ is
then given by
\[
t=\frac1{2\pi i}\log c_2+\frac1{2\pi i}\log z+\frac{g(z)}{f_0(z)}
\]
($\frac1{2\pi i}\log c_2$ is the second ``constant of integration'')
and the associated parameter $q$ is
\[
q=e^{2\pi it}=c_2ze^{g/f_0}.
\]

Let us define
\[ \delta(z)=1+z\frac d{dz}\left(\frac{g(z)}{f_0(z)}\right), \]
so that
\[
\frac{dq}{dz}=c_2\delta(z)e^{g/f_0}.
\]
Then by the chain rule,
\[
\frac{dz}{dt}=\frac{dq/dt}{dq/dz}=\frac{2\pi iz}{\delta(z)}.
\]
It follows that the gauge-fixed value of $\kappa_{ttt}$ is
\[
\kappa_{ttt}=\left(\frac{dz}{dt}\right)^3\kappa_{zzz}
=\frac{c_1}{(\delta(z))^3(z-\lambda)(f_0(z))^2}.
\]

Finally we express this normalized $\kappa_{ttt}$ as a power series
in $q$.  The constants $c_1$ and $c_2$ have
 yet to be determined; however, we
can define
\begin{equation} \label{h0z}
h_0(z)=\frac{1}{(\delta(z))^3(z-\lambda)(f_0(z))^2}
\end{equation}
\begin{equation} \label{hjz}
h_j(z)=\frac1{\delta(z)e^{g/f_0}}\cdot\frac{dh_{j-1}(z)}{dz}
\end{equation}
and find that
\[
h_j(z)=\frac{(c_2)^j}{c_1}\left(\frac d{dq}\right)^j\kappa_{ttt},
\]
so that
\[
\kappa_{ttt}=\sum_{j=0}^{\infty} \frac{c_1}{(c_2)^j}\,
\frac{h_j(0)}{j!}\, q^j.
\]

\begin{proposition}
The numbers $h_j(0)$ are rational numbers.
\end{proposition}

\begin{pf}
The coefficient matrix $A(z)$ in the vector equation \eqref{w0eqn}
has entries in $\Bbb Q(z)$; if written out in power series, all
the power series coefficients will be rational numbers.  Finding
a power series solution to \eqref{w0eqn} then involves solving linear
equations with rational coefficients at each step:  the solutions
will be rational.  Thus, $w_0(z)$ and $f_0(z)$ are power series in
$z$ with rational coefficients.

Similarly, $v(z)$ and $g(z)$ are power series with rational coefficients,
since they come from equation \eqref{veqn}.  Furthermore, since
exponentiating a power series with rational coefficients (whose constant
term is zero) again gives a power series with rational coefficients,
$e^{g/f_0}$ and $\delta(z)$ are power series in $z$ with rational
coefficients.

But then by \eqref{h0z},
$h_0(z)$ is clearly a power series in $z$ with rational
coefficients; similarly for $h_j(z)$ by \eqref{hjz}.  It follows that
each $h_j(0)$ is a rational number.
\end{pf}

\section{Choosing the constants and predicting the
numbers of rational curves}

Calabi-Yau threefolds with $h^{2,1}=1$ are conjectured to be the ``mirrors''
of other Calabi-Yau threefolds with $h^{1,1}=1$.  In the four examples
we have considered, this mirror property can be realized by a construction
of Greene and Plesser \cite{greene-plesser}.  The threefolds ${\cal W}_\psi$
are mirrors of threefolds ${\cal M}\subset\Bbb P^{(k_0,\dots,d_4)}$,
which are hypersurfaces of weighted degree $k=\sum k_j$.  The Picard
group of ${\cal M}$ is cyclic, generated by some ample divisor $H$.

Mirror symmetry predicts that the $q$-expansion of the gauge-fixed
Yukawa coupling
\[
\kappa_{ttt}=a_0+a_1q+a_2q^2+\dotsb
\]
will have integers as coefficients.  Moreover, by a formula conjectured
in \cite{pair} and established in \cite{psa-drm}, if this $q$-expansion
is written in the form
\begin{equation}  \label{formla2}
\kappa_{ttt} = n_0 +
\sum_{j=1}^\infty \frac{n_jj^3q^j}{1-q^j}
= n_0 + n_1q + (2^3n_2 + n_1)q^2 + \dotsb .
\end{equation}
then the coefficients $n_j$ are also integers.  The first term $n_0$
is predicted to coincide with $H^3$ (the absolute degree of ${\cal M}$),
and $n_j$ is predicted to be the number of rational curves $C$ on
${\cal M}$ with $C\cdot H=j$, assuming that all rational curves on ${\cal M}$
are disjoint and have normal bundle ${\cal O}(-1)\oplus{\cal O}(-1)$.

These two predictions can be used to choose the constants of integration
in our examples.
First, the absolute degree $d$ is the lowest order term which appears
in the polynomial $Q(x,\psi)$; to ensure that $n_0=d$ we must take
$c_1=-\lambda d$.  Second, the formula~\eqref{formla2} puts very
strong divisibility constraints on the coefficients $a_j$, and it
seems likely that there will be a unique choice of $c_2$ which
satisfies all of these constraints.

\begin{table}[t]
\renewcommand{\arraystretch}{1.2}
\begin{center}
\begin{tabular}{|c|rrrrr|} \hline
$k$ & $n_0$ & $n_1$ & $n_2$ & $n_3$ & $n_4$ \\ \hline
 5 & 5 & 2875 & 609250 & 317206375 & 242467530000 \\
 6 & 3 & 7884 & 6028452 & 11900417220 & 34600752005688 \\
 8 & 2 & 29504 & 128834912 & 1423720546880 & 23193056024793312 \\
10 & 2 & 462400 & 24431571200 & 3401788732948800 & 700309317702649312000 \\
\hline \end{tabular} \end{center}
\caption{The predicted numbers of curves.}
\label{tab3}
\end{table}

We have calculated the first 20 coefficients (using {\sc mathematica})
in each of our four examples.  There does indeed appear to be a unique
choice for $c_2$ which produces integers for $n_1,\dots n_{20}$:
that choice turns out to be $c_2=k^{-k}$ in each of our examples.
 Making this choice leads to the values for $n_j$
displayed in table~\ref{tab3}.

Table~\ref{tab3} therefore contains predictions about numbers of rational
curves on the weighted projective hypersurfaces.  For a general
hypersurface in ${\cal M}\subset\Bbb P^{(k_0,\dots,d_4)}$ of
degree $k=\sum k_j$,
the prediction is that there should be $n_j$ rational curves $C$ with
$C\cdot H=j$, where $H$ generates $\Pic({\cal M})$.

The first line of
the table reproduces the predictions made by Candelas et al.\ about
quintic threefolds.  Several of these have been verified:  the number
of lines was known classically, the number of conics was computed by
Katz \cite{katz}, and the number of twisted cubics $n_3$ has recently
been computed by Ellingsrud and Str\o mme \cite{ell-str}---all of these
results agree with the predictions.

Of the remaining predictions in the table, we have only checked one.
Each hypersurface from the third family (the case $k=8$) can be regarded as
a double cover of $\Bbb P^3$ branched on a  surface of degree 8.
The entry 29504 in the third line of the table can be interpreted as
follows:  for a general surface of degree 8 in $\Bbb P^3$, there should
be 14752 lines which are 4-times tangent to the surface.  (These lines
will then split into pairs of rational curves on the double cover.)
After we had obtained this number,
Steve Kleiman was kind enough to locate a
$19^{\text{th}}$-century formula of
Schubert \cite[Formula 21, p. 236]{schubert}, which
states that
the number of lines in $\Bbb P^3$ 4-times tangent to a general
surface of degree
$n$ is
\[  \frac1{12}n(n-4)(n-5)(n-6)(n-7)(n^3+6n^2+7n-30)  .
\]
Substituting $n=8$, we find the predicted number 14752.

\section*{Acknowledgements}

Several of the ideas explained in this paper arose in conversations with
Sheldon Katz---it is a pleasure to acknowledge his contribution.  I also
benefited greatly from the chance to interact with physicists which
was provided by the Mirror Symmetry Workshop.  I wish to thank the M.S.R.I.
as well as the organizers and participants in the Workshop.

This work was partially supported by NSF grant DMS-9103827.



\begin{thebibliography}{10}

\bibitem{psa-drm}
P.~S. Aspinwall and D.~R. Morrison, {\em Topological field theory and rational
  curves}, Oxford and Duke preprint OUTP-91-32P, DUK-M-91-12, October 1991.

\bibitem{blok-var}
B.~Blok and A.~Varchenko, {\em Topological conformal field theories and the
  flat coordinates}, IAS preprint IASSNS-HEP-91/5, January 1991.

\bibitem{buchberger}
B.~Buchberger, {\em {G}r\"obner bases: An algorithmic method in polynomial
  ideal theory}, Multidimensional Systems Theory
  (N.~K. Bose, ed.), D. Reidel, Dordrecht, Boston, Lancaster, 1985,
pp.~184--232.

\bibitem{cad-ferr}
A.~C. Cadavid and S.~Ferrara, {\em {P}icard-{F}uchs equations and the moduli
  space of superconformal field theories}, Phys. Lett. B {\bf 267} (1991),
  193--199.

\bibitem{pair}
P.~Candelas, X.~C. de~la Ossa, P.~S. Green, and L.~Parkes, {\em A pair of
  {C}alabi-{Y}au manifolds as an exactly soluble superconformal theory},
Phys. Lett. B {\bf 258} (1991), 118--126;
  Nuclear Phys. B {\bf 359} (1991), 21--74.

\bibitem{cls}
P.~Candelas, M.~Lynker, and R.~Schimmrigk, {\em {C}alabi-{Y}au manifolds in
  weighted {$\Bbb P_4$}}, Nuclear Phys. B {\bf 341} (1990), 383--402.

\bibitem{codlev}
E.~A. Coddington and N.~Levinson, {\em Theory of ordinary differential
  equations}, McGraw-Hill, New York, Toronto, London, 1955.

\bibitem{deligne}
P.~Deligne, {\em Equations diff\'erentielles \`a points singuliers
  r\'eguliers}, Lecture Notes in Math., vol. 163, Springer-Verlag, Berlin,
  Heidelberg, New York, 1970.

\bibitem{dolg}
I.~Dolgachev, {\em Weighted projective varieties}, Group Actions and Vector
  Fields  (J.~B. Carrell, ed.), Lecture Notes in
  Math., vol. 956, Springer-Verlag, Berlin, Heidelberg, New York, 1982,
pp.~34--71.

\bibitem{dwork}
B.~Dwork, {\em On the {Z}eta function of a hypersurface, {II}}, Ann. of Math.
  (2) {\bf 80} (1964), 227--299.

\bibitem{ell-str}
G.~Ellingsrud and S.~A.~Str\o mme, {\em The number of twisted cubic curves on
  the general quintic threefold}, preprint, 1991.

\bibitem{ferrara}
S.~Ferrara, {\em {C}alabi-{Y}au moduli space, special geometry and mirror
  symmetry}, Modern Phys. Lett. A {\bf 6} (1991), 2175--2180.

\bibitem{greene-plesser}
B.~R. Greene and M.~R. Plesser, {\em Duality in {C}alabi-{Y}au moduli space},
  Nuclear Phys. B {\bf 338} (1990), 15--37.

\bibitem{griffiths}
P.~A. Griffiths, {\em On the periods of certain rational integrals, {I}}, Ann.
  of Math. (2) {\bf 90} (1969), 460--495.

\bibitem{gr-bull}
\bsame, {\em Periods of integrals on algebraic manifolds: Summary of main
  results and discussion of open problems}, Bull. Amer. Math. Soc. {\bf 76}
  (1970), 228--296.

\bibitem{transcendental}
\bsame, ed., {\em Topics in transcendental algebraic geometry}, Ann. of
  Math. Stud., vol. 106, Princeton University Press, Princeton, 1984.

\bibitem{nkatz}
N.~M. Katz, {\em Nilpotent connections and the monodromy theorem: Applications
  of a result of {T}urrittin}, Inst. Hautes {\'E}tudes Sci. Publ. Math. {\bf
  39} (1970), 175--232.

\bibitem{katz}
S.~Katz, {\em On the finiteness of rational curves on quintic threefolds},
  Compositio Math. {\bf 60} (1986), 151--162.

\bibitem{lsw}
W.~Lerche, D.-J. Smit, and N.~P. Warner, {\em Differential equations for
  periods and flat coordinates in two dimensional topological matter theories},
  preprint LBL-31104, UCB-PTH-91/39, USC-91/022, CALT-68-1738, July 1991.

\bibitem{markushevich}
D.~G. Markushevich, {\em Resolution of singularities (toric method)}, appendix
  to: D. G. Markushevich, M. A. Olshanetsky, and A. M. Perelomov, {\em
  Description of a class of superstring compactifications related to
  semi-simple {L}ie algebras}, Comm. Math. Phys. {\bf 111} (1987), 247--274.

\bibitem{guide}
D.~R. Morrison, {\em Mirror symmetry and rational curves on quintic threefolds:
  A guide for mathematicians}, Duke preprint DUK-M-91-01, July 1991.

\bibitem{roan1}
S.-S. Roan, {\em On the generalization of {K}ummer surfaces}, J. Differential
  Geom. {\bf 30} (1989), 523--537.

\bibitem{roan2}
\bsame, {\em On {C}alabi-{Y}au orbifolds in weighted projective spaces},
  Internat. J. Math. {\bf 1} (1990), 211--232.

\bibitem{roanpf}
\bsame, {\em The mirror of {C}alabi-{Y}au orbifold}, Max-Planck-Institut
  preprint MPI/91-1, to appear in Internat. J. Math., 1991.

\bibitem{schubert}
H.~C.~H. Schubert, {\em {K}alk\"ul der abz\"ahlenden {G}eometrie}, 1879,
  reprinted with an introduction by S. Kleiman, Springer-Verlag, 1979.

\bibitem{steen}
J.~Steenbrink, {\em Intersection form for quasi-homogeneous singularities},
  Compositio Math. {\bf 34} (1977), 211--223.

\bibitem{yau}
S.~T. Yau, {\em On {C}alabi's conjecture and some new results in algebraic
  geometry}, Proc. Nat. Acad. Sci. U.S.A. {\bf 74} (1977), 1798--1799.

\bibitem{yu}
Y.~Yu, {\em An improvement on the {G}r\"obner basis algorithm}, in preparation.

\end{thebibliography}

\makeatletter \renewcommand{\@biblabel}[1]{\hfill#1.}\makeatother
\newcommand{\bsame}{\leavevmode\hbox to3em{\hrulefill}\,}

\trailer

\end{document}